\documentclass{ws-jml-fixed}

\usepackage{amsmath,amssymb,amsfonts}
\usepackage{amsthm}
\usepackage{xspace}
\usepackage{booktabs}
\usepackage{color,epsfig}
\usepackage{burdges-ws-jml-fixed}
\usepackage{calc}

\newcommand\EE{\mathcal{E}}
\newcommand\oo{\hbox{$^\circ$}\xspace}

\newcommand\Qq{\mathbb Q}
\newcommand\Zz{\mathbb Z}

\renewcommand\Submitted[1][]{Submitted}

\begin{document}

\markboth{Jeffrey Burdges, Gregory Cherlin}{A generation theorem for groups of finite Morley rank}

\def\version{August 2007}
\title{A generation theorem for \\ groups of finite Morley rank}

\author{Jeffrey Burdges\thanks{Supported by NSF postdoctoral fellowship DMS-0503036, a Bourse Chateaubriand postdoctoral fellowship, and DFG grant Te 242/3-1.}}

\address{Institut Camille Jordan, Universit\'e de Lyon-I, Lyon, France}

\author{Gregory Cherlin\thanks{Supported by NSF grant DMS-0100794 and DMS-0600940}}

\address{Mathematics Department, Rutgers University, Piscataway, NJ, USA}

\maketitle

\begin{abstract}
We deal with two forms of the ``uniqueness cases'' in the classification of
large simple $K^*$-groups of finite Morley rank of odd type, where large means
the $m_2(G)$ at least three.  This substantially extends results known for
even larger groups having \Prufer 2-rank at least three, to cover the two
groups $\PSp_4$ and $\G_2$.  With an eye towards distant developments,
we carry out this analysis for $L^*$-groups which is substantially broader
than the $K^*$ setting.
\end{abstract}


\section*{Introduction}

\subsection*{Groups of finite Morley rank}

A group of finite Morley rank is a group equipped with a notion
of dimension satisfying various natural axioms \cite[p.~57]{BN};
These groups arise naturally in model theory,
 expecially geometrical stability theory.  
The main examples are algebraic groups over algebraically closed
fields, where the notion of dimension is the usual one,
as well as certain groups arising in applications of model theory to
diophantine problems, where the notion of dimension
comes from differential algebra rather than algebraic geometry.

The structural analysis of groups of finite Morley rank is dominated
by the Algebraicity Conjecture (Cherlin/Zilber), which states that
a simple group of finite Morley rank should, in fact, be a Chevalley
group over an algerbaically closed field, i.e., an algebraic group. 
This is a strong conjecture, which asserts that the classification of the
simple algebraic groups can be carried out using only their coarsest
properties. 

It is known that a counterexample to the Algebraicity Conjecture 
containing at least one involution (i.e., element of order two)
must have {\em odd type}.
This means that the connected component of a Sylow 2-subgroup
 is a 2-torus (i.e., a divisible abelian group).
This reduction to odd type is the result of a large body of work
 presented in \cite{ABC_EvenType}. 
An equivalent condition is that the group $G$ has
 infinite Sylow 2-subgroups of finite {\em 2-rank}.
Here one defines the 2-rank $m_2(G)$ to be the supremum
 of the ranks of elementary abelian 2-subgroups of $G$.

\subsection*{Generic vs.~thin}
Our aim in the present work is to lay the foundations of the theory of
``generic'' simple groups of finite Morley rank as broadly as possible,
taking all groups with $m_2(G) \ge 3$ as generic.
A more concrete goal, not reached here, is the following.

\begin{namedtheorem}{Generic Algebraicity Conjecture}
Let $G$ be a simple $K^*$-group of finite Morley rank 
with $m_2(G) \geq 3$.  Then $G$ is a Chevalley group over an
algebraically closed field, of Lie rank at least two.
\end{namedtheorem}

A $K^*$-group is one whose proper definable simple sections are
algebraic. A minimal nonalgebraic simple group of finite Morley rank
would be a $K^*$-group.

The point of stating this conjecture separately under the hypothesis
$m_2(G)\ge 3$ is that this is the natural target for methods of a
{\em general} character. 
Some results of this type are already known, but under more restrictive
assumptions. 
Building on this work, we aim to broaden the class of groups
which may be treated as generic to its widest natural limit, 
leaving over as little as possible for separate
consideration on a case by case basis. 

The philosophy in the study of finite simple groups, which has been
transferred also to the case of simple groups of finite Morley rank,
is that large groups are those which can be handled by reasonably
general methods, and small groups are those which require close
analysis by ad hoc methods. The terms in use in finite group theory
for ``large'' and ``small'' are ``generic'' and ``thin'', 
with an intermediate category, ``quasithin''. 
The precise definitions even in finite group
theory remain in flux; they are chosen to match the details of
specific proof strategies, and in the two generations of the proof of
the classification of finite simple groups there has been considerable
variation in the treatment of borderline cases. 
But tinkering aside, these are robust
notions: generic groups are for the most part groups which turn out to
be Chevalley groups of Lie rank at least 3, quasithin more or less Lie
rank 2, and the thin groups are populated by Lie rank 1 and more
pathological groups which offer comparably little internal structure. 

Our goal here is to treat a group $G$ as {\em generic}
 if $m_2(G)\ge 3$; thus we treat some quasithin cases.
Stricter notions of genericity are defined by considering the
structure of the Sylow 2-subgroup $S$, whose {\em normal} 2-rank
$n_2(S)$ and {\em Pr\"ufer} 2-rank $\pr_2(S)$ are defined by taking
the maximal rank of a normal elementary abelian 2-subgroup, or an
elementary abelian 2-subgroup contained in a divisible abelian
2-subgroup, respectively. One sets $n_2(G)=n_2(S)$, $\pr_2(G)=\pr_2(S)$,
which by conjugacy of Sylow 2-subgroups is well-defined.
It is easy to see that $$m_2(G)\ge n_2(G)\ge \pr_2(G)$$

\subsection*{Generic groups and Uniqueness Cases}

The following approximation to the Generic Algebraicity
Conjecture is proved in \cite{Bu07a};
an earlier version, assuming the absence of ``bad
fields'', is due to Borovik.

\begin{namedtheorem}{High Pr\"ufer rank Case}
A simple $K^*$-group of finite Morley rank with
Pr\"ufer 2-rank at least three is algebraic.
\end{namedtheorem}

Pushing the same result through in the case of $m_2(G)\ge 3$ 
involves substantial technical complications.
The general thrust of the argument in \cite{Bu07a} is to first exclude
a ``uniqueness case'' involving the {\em 2-generated core}, and
then to apply systematic methods---signalizer functor theory and 
the Curtis-Tits theorem, primarily---to carry through the necessary
structural analysis and recognition process.

Broadly speaking, a {\em uniqueness subgroup} of a simple group $G$ is
a very large definable subgroup which does not interact much with
other subgroups of $G$. A typical example is the
Borel subgroup of $\PSL_2$. In generic type groups one aims to show
that there can be no proper uniqueness subgroup in one or another
sense. This then becomes the point of departure for a 
systematic study of interactions between subgroups, and eventually
for the recognition of the group, or of a geometry on which the group acts.

In our context there are four notions of uniqueness subgroup that come
into play: 
\begin{enumerate}
\item The groups $\Gamma_V$ where $V$ is elementary abelian $p$-subgroup
  are $$\gen{C\oo(E): E\leq E, [E,V]=p}$$
\item The weak 2-generated core, denoted $\Gamma_{S,2}^0$:
  for $S\le G$ a Sylow 2-subgroup of $G$, this is the definable hull
  of the group generated by all subgroups $$ N(A) $$
  where $A$ varies over rank 2 elementary abelian subgroups of $S$ which
  are contained in rank 3 elementary abelian subgroups of $S$.
\item The 2-generated core, denoted $\Gamma_{S,2}$, is defined
  similarly, but taking into account all elementary abelian subgroups
  $A$ of $S$ of rank 2, including those that are maximal in $S$.
\item Strongly embedded subgroups: a proper definable
  subgroup $M$ of $G$ is strongly embedded if $M$ contains an
  involution while $M\intersect M^g$ contains no involution for
  $g\notin M$.
\end{enumerate}

In the definition of the weak 2-generated core, the superscript $0$ is
not a reference to the connectedness of any group, but rather an
allusion to a connected component in an associated graph of
elementary abelian $p$-subgroups. This graph is introduced explicitly
in \S\ref{sec:Generation}.

The results we aim at are the following. The second turns out to be a
corollary of the first.

\begin{namedtheorem}{Generation Theorem I}
Let $G$ be a simple $K^*$-group of finite Morley rank and odd type, 
with $m_2(G) \geq 3$. Then $\Gamma_V = G$
whenever $V$ is an elementary abelian subgroup of rank 2.
\end{namedtheorem}

\begin{namedtheorem}{2-Generated Core Theorem I}
Let $G$ be a simple $K^*$-group of finite Morley rank and odd type, 
with $m_2(G) \geq 3$. Then
$\Gamma_{S,2}^0 = G$ for any Sylow 2-subgroup $S$ of $G$.
\end{namedtheorem}

These results will be given in a more general form as Theorem
\ref{K*-generation}, and in the forms just stated as Corollary
\ref{GenerationI}.

The first of these results is a generation theorem: $G$ is generated
by the connected components of centralizers $C\oo(v)$ 
when $v$ varies over involutions in an elementary
abelian $2$-subgroup of order $4$. This points the way toward
the identification of the group $G$ in terms of the structure of
centralizers of involutions.
The meaning of the second result is less transparent,
 but using signalizer functor theory it can be shown that 
a proper weak 2-generated core $\Gamma^0_{S,2}$ arises whenever
the centralizer $C(i)$ of in involution is ``not sufficiently reductive''
(meaning that $O(C(i))$ is ``large'' enough to be ``unipotent'').
So for generic $G$, these statements give both a generation property
by related centralizers of involutions, and 
a weak reductivity condition applying to those centralizers. 
Combined, in Pr\"ufer rank at least three these results lead eventually to
generation of $G$ by quasisimple components of centralizers of
involutions, and then to a full description of $G$. Here the
Curtis-Tits theorem simplifies matters considerably; in Pr\"ufer rank
two other methods come into play (one may consider this a quasithin
case: generic arguments carry one some distance, then one must enter
into the consideration of a certain number of special configurations).

In the proofs of our results we follow the line of 
\cite{BBN04}, showing that our hypothetical counterexample $G$
is a minimal simple group of Pr\"ufer 2-rank at least two, with a strongly
embedded subgroup.  Since such a configuration is impossible by
\cite{BCJ}, this suffices. The line of argument goes as follows: if
$\Gamma_V$ is the offending uniqueness subgroup (namely, $\Gamma_V<G$)
then $N(\Gamma_V)$ turns out to be strongly embedded; this is the
archetypal uniqueness condition. Further analysis then shows that $G$
is a minimal simple group, reaching a contradiction to \cite{BCJ}.

A technical point of considerable importance is some new
information about ``semisimple torsion'' (e.g., $2$-elements
in odd type groups) which has become available only recently. 
This development begins with \cite{BBC} and has been carried further
in \cite{BC08a}, with particular attention to results which are
useful in the present analysis.  We are able to put aside a number
of pathological configurations on this basis, greatly simplifying
the flow of the argument.

\subsection*{Variation: $p$}

We will take some pains to place our two main theorems in a substantially
broader context. In the first place, we will replace the prime 2 by a
general prime $p$.  For this, we ask that $G$ have $p^\perp$ type meaning
no unipotent $p$-subgroup, i.e.\ no connected definable $p$-subgroup of
bounded exponent.  This requires some additional work, some of it already
carried out in \cite{BC08a}.

Our two uniqueness case theorems take on the following forms in general.

\begin{namedtheorem}{Theorem \ref{K*-generation}}
Let $G$ be a simple $K^*$-group of finite Morley rank and $p^\perp$ type, 
with $m_p(G) \geq 3$ and $\pr_2(G)\ge 2$. 
Then $\Gamma_V = G$ for any elementary abelian $p$-subgroup of rank two.
In particular, $\Gamma_{S,2}^0 = G$ for any Sylow $p$-subgroup $S$ of
$G$. 
\end{namedtheorem}

Two important technical points which come up here are the following.

\begin{namedtheorem}{Theorem \ref{Prank2}}
Let $G$ be a connected group of finite Morley rank and $p^\perp$ type
with $m_p(G) \geq 3$.  Then $\pr_p(G)\ge 2$.
\end{namedtheorem}

This is particularly important for $p=2$.

\begin{namedtheorem}{Theorem \ref{max_el_abelian}}
Let $G$ be a connected group of finite Morley rank and $p^\perp$ type.
with $m_p(G) \geq 3$.  Then
$m_p(V)\geq 3$ for any maximal elementary abelian $p$-subgroup $V \leq G$.
\end{namedtheorem}

Theorem \ref{Prank2} should always be borne in mind.
We would not relish being forced to treat cases of Pr\"ufer rank 1
within the ``generic'' framework. 

\subsection*{Variation: $L^*$}

Returning to the case $p=2$, we would also like to dispense with the
restriction to $K^*$-groups. If all simple $K^*$-groups are algebraic,
then the same follows for all groups of finite Morley rank. But we
would prefer to deal with groups of odd type in a way which does not
require a complete prior analysis of groups of degenerate type.
Indeed, currently we do not even know if simple algebraic groups of
finite Morley rank 
are necessarily $K$-groups. 

However, non-algebraic definable subgroups of such algebraic groups
of finite Morley rank necessarily have degenerate type \cite{BB_Linear}.
So the natural class of groups to consider, called {\em $L^*$-groups},
are those whose proper definable connected simple sections are either
algebraic or of degenerate type, the latter without involutions by \cite{BBC}.
The analysis of $L^*$-groups then presupposes some more limited prior
analysis of $L$-groups, where an $L$-group satisfies the same condition
for {\em all} definable connected simple sections.

Versions of our main results will be proved under the $L^*$
hypothesis in a form that then yields the stated results immediately
under the stronger $K^*$-hypothesis.  We give this with an eye to
the future---we are a long way from having the same kind of theory
in an $L^*$ context that we have for $K^*$-groups, but at least the
present chapter of the theory has a satisfactory development at that
level of generality.

\begin{namedtheorem}{Theorem \ref{thm:strongembedding}}
Let $G$ be a simple $L^*$-group of finite Morley rank 
 and odd type,  with $m_2(G) \geq 3$.
Suppose that 
$$\hbox{$\Gamma_V<G$ for some $V \leq G$ with $m_2(V) \geq 2$}$$
Then $N(\Gamma_V)$ is a strongly embedded subgroup of $G$.
Furthermore $N(\Gamma_V)=\Gamma_{S,2}(G)$ and $\Gamma_V=\Gamma_{S,2}\oo$.
\end{namedtheorem}

Actually, the statement of Theorem \ref{thm:strongembedding} below is
slightly more general. We replace the $L^*$ condition by an
assumption that a generation statement applies to proper connected
sections of $G$.
The corresponding generation statement for $L$-groups, which is needed
to justify the result as phrased above, is found in Theorem
\ref{generationA}.  

This result may be elaborated on further.
We will call a group of finite Morley rank a {\em $D$-group} if all its
definable connected simple sections are of degenerate type, and a {\em
  $D^*$-group} if the same applies to all proper definable connected
simple sections.
In the context of $K^*$ groups of odd type, a $D^*$ group would be a
minimal connected simple group. 

\begin{namedtheorem}{Theorem \ref{L*-generation}}
Let $G$ be a simple $L^*$-group of finite Morley rank 
 and odd type,  with $m_2(G) \geq 3$.
Let $V$ be an elementary abelian $2$-group of $2$-rank two with
$\Gamma_V<G$. Then
\begin{enumerate}
\item $G$ is a $D^*$-group.
\item $M := N(\Gamma_V)$ is a strongly embedded subgroup of $G$.
\item The Sylow 2-subgroups of $G$ are connected.
\item If $r$ is the least prime divisor of $|M/M\oo|$, 
then $G$ contains a nontrivial unipotent $q$-subgroup for some $q\le r$.
\end{enumerate}
\end{namedtheorem}

In this formulation, the result as stated is very close to the version
given in \cite{BBN04} under more restrictive conditions.
In the $K^*$ case one 
reaches a contradiction by applying \cite{BCJ},
which says that minimal connected simple groups of finite Morley rank
and odd type with $\pr_2(G) \geq 2$ never have strongly embedded subgroups.
The corresponding result is not known in the $D^*$ context, and should
be difficult.

Finally, to close the gap between $\Gamma_V$ and $\Gamma_{S,2}$, one
more result is required. 

\begin{namedtheorem}{Corollary \ref{cor:2-gencore}}
Let $G$ be a simple $L^*$-group of finite Morley rank 
 and odd type,  with $m_2(G) \geq 3$.
Then $\Gamma_V \le \Gamma_{S,2}$.
\end{namedtheorem}

Again, the statement actually given as Corollary \ref{cor:2-gencore}
replaces the $L^*$-group hypothesis by an inductive condition on
sections.

We have the following ``characteristic zero'' corollary.

\begin{namedtheorem}{Corollary}
Let $G$ be a simple $L^*$-group of finite Morley rank with $m_2(G) \geq 3$,
which has no unipotent torsion.  Then
\begin{enumerate}
\item $\Gamma_V = G$ whenever $m_2(V) \geq 2$; equivalently,
\item $\Gamma_{S,2}^0 = G$ for any Sylow 2-subgroup $S$ of $G$.
\end{enumerate}
\end{namedtheorem}

At this stage it is natural to focus on the $K^*$ case, and to
ignore the further subtleties that intervene at a more general level,
but even in that case something very much like our $L^*$ version must
be proved on the way toward a contradiction via \cite{BCJ}.

\subsection*{Contents}

We begin by proving Theorems \ref{Prank2} and \ref{max_el_abelian}.
These come quickly out of the general theory developed in \cite{BBC}, 
We also give
some $L$-group theory which is parallel to what is known in the
$K$-group case, in preparation for the proof 
of our generation theorems.

After that we give a preparatory generation theorem for
$L$-groups when $p=2$, and for $K$-groups when $p\ne 2$.
This takes considerable analysis, involving a reduction to the simple
algebraic case.
The statement runs as follows. 

\begin{namedtheorem}{Theorem \ref{generationA}}
Fix a prime $p$.
Let $H$ be a connected group of finite Morley rank.
If $p=2$ suppose that $H$ is an $L$-group of $2^\perp$ type;
if $p>2$ suppose that $H$ is a $K$-group of $p^\perp$ type.
Let $E$ be a finite nontrivial elementary abelian $p$-group acting definably on $H$.
Then $$ \Gamma_V(H) = H $$.
\end{namedtheorem}

After that we turn to the generation theorem for simple $K^*$-groups
and $L^*$-groups.
We also show that the corresponding results for the weak 2-generated
core $\Gamma_{S,2}^0$ follow from the $\Gamma_V$ versions.

After these preparations, our next result is strong embedding for
$N(\Gamma_V)$ when $\Gamma_V<G$.  In section \S\ref{sec:Altseimer} we prove
that a strongly embedded subgroup of our group $G$ is then a $D$-group.
The argument is based in part on one given by Christine Altseimer in 
\cite{AltseimerStrEmb} in the $K^*$ context.

\begin{namedtheorem}{Theorem \ref{thm:Altseimer}}
Let $G$ be a simple $L^*$-group of finite Morley rank and odd type.
Then a strongly embedded subgroup $M$ of $G$ is
 either a $D$-group, or else satisfies $M\oo/\hat O(M\oo) \iso \SL_2$.
\end{namedtheorem}  

In the last section we complete our analysis by combining
Theorems \ref{thm:strongembedding} and \ref{thm:Altseimer},
finally proving that $G$ is a $D^*$-group.
Then \cite{BCJ} provides a contradiction in the $K^*$ case, proving
the generation theorem in that case.



\section{Preliminaries}

We will prove the following two results in the present section.

\begin{theorem}\label{Prank2}
Let $G$ be a connected group of finite Morley rank and $p^\perp$ type
with $m_p(G) \ge 3$.
If $p>2$, then $\pr_p(G) \ge 3$. If $p=2$, then $\pr_2(G) \ge 2$.
\end{theorem}

\begin{theorem}\label{max_el_abelian}
Let $G$ be a connected group of finite Morley rank and $p^\perp$ type
with $m_p(G) \geq 3$. Then any maximal elementary abelian $p$-subgroup
 $V \leq G$ has $p$-rank at least 3.
\end{theorem}

We will also develop some general $L$-group theory which comes into
play in our analysis.

\subsection{Semisimple torsion}

We will rely on results on $p$-torsion in groups of $p^\perp$ type
developed in \cite{BC08a}, notably the following.

\begin{fact}[{\cite{BC08a}}]\label{torality}
Let $G$ be a connected group of finite Morley rank of $p^\perp$ type.
Let $a\in G$ be a $p$-element, and $T\le C(a)$ a maximal $p$-torus in $C(a)$. 
Then $a\in T$.
\end{fact}

This is commonly applied in two different ways. First and foremost, 
any $p$-element will be {\em toral} (lie in a $p$-torus). 
Secondly, any $p$-element in the centralizer of a maximal $p$-torus $T$
of $G$ will lie in $T$, or in other words any $p$-element in $N(T)\setminus T$
will act nontrivially on $T$ (in the latter formulation, this is
 a statement about the {\em Weyl group} $N(T)/C\oo(T)$). 

We have known that Sylow 2-subgroups exist and are conjugate for sometime \cite{BP}.
However, we have established the conjugacy theorem for Sylow $p$-subgroups
with $p$ odd only recently, and only in the $p^\perp$ case.
By definition, a Sylow $p$-subgroup of a group of finite Morley rank will be
a maximal solvable $p$-subgroup (equivalently, a maximal locally finite
$p$-subgroup---which makes existence clearer).

\begin{fact}[{\cite[Theorem 4]{BC08a}}]
Let $G$ be a group of finite Morley rank of $p^\perp$ type. 
Then Sylow $p$-subgroups of $G$ are conjugate.
\end{fact}

For $p=2$ this is known for arbitrary groups of finite Morley rank.

We will also require some elementary facts about automorphisms of $p$-tori.

\begin{fact}[{\cite[\S3.3]{Be01,BB04}}]\label{endptor}
Let $T$ be a $p$-torus of Pr\"ufer $p$-rank $d$.
Then $\End(T)$ can be faithfully represented as the ring $M_d(\Zz_p)$
of $d \times d$ matrices over the $p$-adic integers $\Zz_p$.
\end{fact}

Underlying this fact is a duality between $T$ and its
so-called {\em Tate module}, which is a free $\Zz_p$-module of rank
$d$ naturally associated with $T$. 
We remark that the action of $M_d(\Zz_p)$ on $T$ is natural
 if $T$ is represented as $C_p^d$ with $C_p$ quasicylic
 (the direct limit of cyclic groups $\Zz/p^n\Zz$);
multiplication of elements of $C_p$ by elements of $\Zz_p$ is well defined,
 since $\Zz_p/p^n\Zz_p=\Zz/p^n\Zz$ for any $n$. 

\begin{lemma}\label{centinvT}
Let $T$ be a nontrivial 2-torus and let $i$ be an involution acting on $T$.
Then either $i$ inverts $T$ or $C_T(i)$ is infinite.
\end{lemma}

\begin{proof}
Viewing $i$ as given by a matrix $E$ in $M_d(\Qq_p)$, if $E-I$ is
noninvertible then $E-I$ has a kernel in its natural action on $\Zz_p^d$.
Fixing a nonzero vector $v\in \Zz_p^d$ annihilated by $E-I$ and working
modulo $p^n$ for large $n$, we obtain arbitrarily large cyclic
subgroups of $T$ centralized by $i$. So in this case $C_T(i)$ is infinite.

On the  other hand if $E-I$ is invertible then since $(E-I)(E+I)=0$
 we have $E=-I$ and $i$ acts by inversion.
\end{proof}

This proof can be given just as easily in the language of finite group theory,
but as we have a description of $\End(T)$ available we prefer to use it.

\begin{lemma}\label{Omega_1(T)}
Let $G$ be a connected group of finite Morley rank, $p$ an odd prime,
$T$ a maximal $p$-torus of $G$, and $a\in N(T)$ a $p$-element
centralizing $\Omega_1(T)$. Then $a$ centralizes $T$.
\end{lemma}

\begin{proof}
Suppose that $a\notin C(T)$. Passing to a power of $a$ we may suppose
that $a^p$ centralizes $T$. So $a$ induces an automorphism of $T$ of
order $p$. Now the action of $a$ on $T$ is represented by a matrix $M$ 
over $\Zz_p$ which acts trivially on $\Omega_1(T)$ and therefore is
congruent to the identity $I$ modulo $p$. If $M \neq I$, write
$$M=I+p^k A$$
with $k$ maximal and $A$ a matrix over $\Zz_p$. Now $M^p=I$, so as
$p>2$ we have by the binomial theorem
$$I\equiv I+p^{k+1}A \mod\ p^{2k+1}$$
This forces $A$ to be divisible by $p$, a contradiction. 
\end{proof}

\begin{lemma}\label{p-aut}
Let $T$ be a $p$-torus of Pr\"ufer rank 1 (a quasicyclic $p$-group),
with $p>2$.
Then $\Aut(T)$ contains no element of order $p$.
\end{lemma}
\begin{proof}
Such an element would correspond to a primitive $p$-th root of unity
in $\Zz_p$, and by the usual argument the cyclotomic polynomial of
degree $p-1$ is irreducible over $\Zz_p$, forcing $p-1\le 1$.
\end{proof}

\begin{fact}\label{Z_meet_o}
Let $P$ be an infinite solvable $p$-subgroup of a group $G$ of finite
Morley rank and $Q$ a normal subgroup of $P$. Then $Q\intersect Z(P)>1$.
\end{fact}

\begin{proof}
The group $P$ is nilpotent by finite \cite{BN}, and the claim is
easily proved for nilpotent by finite $p$-groups. 
\end{proof}

\subsection{Theorem \ref{Prank2}}

We treat the two claims of Theorem \ref{Prank2} separately.

\begin{lemma}
Let $G$ be a connected group of finite Morley rank and $p^\perp$ type,
with $m_p(G)\ge 3$ and $p>2$.  Then $\pr_p(G)\ge 3$.
\end{lemma}

\begin{proof}
We let $T$ be a maximal $p$-torus of $G$ and $S$ a Sylow $p$-subgroup
of $G$ containing $T$. By the conjugacy of Sylow $p$-subgroups of $G$, 
$S$ contains an elementary  abelian $p$-subgroup $E$ with $m_p(E) =
3$. We choose $E$ so as to maximize $A=E\intersect \Omega_1(T)$.
Suppose toward a contradiction that $m_p(T)\le 2$.

Suppose first that $A<\Omega_1(T)$. Fix $t\in \Omega_1(T)\setminus
A$ such that the image of $t$ in $(\Omega_1(T)E)/A$ is in the center
of this group.
Looking at the commutator map $[t,\cdot]:E\to A$, since $A$ has
$p$-rank at most $1$ we find $m_p(C_E(t))\ge 2$.
Thus we may replace $E$ by a subgroup of $p$-rank $3$
containing $A$ and $t$,
contradicting the maximality of $A$.

So $A=\Omega_1(T)$. Hence $E$ acts trivially on $\Omega_1(T)$, and 
by Lemma \ref{Omega_1(T)} $E$ centralizes $T$. Then by
Fact \ref{torality} we have $E\le T$, and our claim  follows.
\end{proof}

\begin{lemma}
Let $G$ be a connected group of finite Morley rank and odd type, with 
$m_2(G)\ge 3$.
Then $\pr_2(G)\ge 2$.
\end{lemma}

\begin{proof}
As in the previous argument we begin with $T$ a maximal $2$-torus of
$G$ and $E\le N(T)$ an elementary abelian $2$-subgroup of $2$-rank
$3$, chosen so as to maximize $A=E\intersect \Omega_1(T)$. 

If $T=1$ then we contradict Fact \ref{torality} (as well as older
results of \cite{BBC} on which the proof of this fact depends), 
so we need only consider the case in which 
$$m_2(T)=1$$

Now each element of $E$ either centralizes or inverts $T$ by Lemma
\ref{centinvT}. So $|C_E(T)|\ge 4$. On the other hand $C_E(T)\le T$ by
Fact \ref{torality}, so $|C_E(T)|\le 2$, a contradiction.
\end{proof}

So Theorem \ref{Prank2} is proved.

\subsection{Theorem \ref{max_el_abelian}}

\begin{proof}[Proof of Theorem \ref{max_el_abelian}]
We have a connected group $G$ of finite Morley rank and $p^\perp$ type,
with $m_p(G) \geq 3$. We claim that any elementary abelian
$p$-subgroup is contained in one of $p$-rank at least three.
We know that $\pr_p(G)\ge 2$ by Theorem \ref{Prank2}.

Suppose toward a contradiction that some maximal elementary abelian
$p$-subgroup $V$ of $G$ has $m_p(V) \leq 2$.
Let $S$ be a Sylow $p$-subgroup of $G$ containing $V$.
Let $T=S\oo$.
By the conjugacy of Sylow $p$-subgroups, we have $m_p(S)\ge 3$.

Let $A=\Omega_1(Z(S))$.
By maximality of $V$, we have $A\le V$. If
$m_p(A)\ge 2$ we then have $V=A\le Z(S)$, and as $m_p(S)\ge
3$ we reach a contradiction to the maximality of $V$.
So $A$ is cyclic of order $p$.
As $\Omega_1(T)\intersect Z(S)>1$ by Fact \ref{Z_meet_o},
 we have $\Omega_1(T)\intersect Z(S)=A\le V$.

If $V\le T$, then by maximality of $V$ we have
$V = \Omega_1(T)\normal S$.
Choose $E\leq S$ elementary abelian with $p$-rank 3.
Then $[E,V]\le A$ and viewing this as a bilinear map from $E\times V$
to $A$, we find that $C_E(V)$ has $p$-rank at least 2.
By maximality of $V$ we have $C_E(V) \leq V$ and thus $V = C_E(V) \leq E$,
 which again contradicts the maximality of $V$. 

So we find 
$$V\not\le T\leqno(*)$$
Choose $v \in V \setminus T$.
Then $V = \gen{A,v}$.
If $m_p(T)\ge 3$ then looking at commutation with $v$ as a map from
$\Omega_1(T)/A$ to $A$, we find an element $t\in \Omega_1(T)\setminus
A$ centralizing $V$ and contradict the maximality of $A$. So
$m_p(T)=2$ and in particular $p=2$ by Theorem \ref{Prank2}.

By maximality of $V$,
$v$ does not centralize $\Omega_1(T)$, and hence
$v$ does not invert $T$.
By Lemma \ref{centinvT}, $T_0 := C_T\oo(v)$ is nontrivial.
Note that $A=\Omega_1(T_0)$.
Let $T_1$ be a maximal torus of $C(v)$ containing $T_0$.
By Fact \ref{torality},
 $v\in T_1$. So $V=\gen{A,v}\le T_1$.

We have shown that $V$ lies inside a $p$-torus of $G$.
But then we could take $T$ to be such a $p$-torus,
contradicting $(*)$.
\end{proof}

\subsection{Structure of $L$-groups}
\newcommand\alg{{\mbox{\rm \small alg}}}

In this section, we give some variations on \cite[Thm.~5.12]{Bo95},
adapted to an $L$-group context. 
We use the notation $\hat O(H)$ for the largest definable connected
normal subgroup of $H$ without involutions, $E(H)$ for the product
of the quasisimple subnormal subgroups $L$ of $H$,
 $F(H)$ for the Fitting subgroup, and
 $F^*(H)=F(H)E(H)$ for the generalized Fitting subgroup.
By \cite[Lemmas 7.9, 7.10, 7.13]{BN},
the product $F(H)E(H)$ is a central product, the group $E(H)$ is a
central product of finitely many normal quasisimple factors, each
definable in $G$, and 
$$C_H(F^*(H))\le F^*(H)$$

We also recall that, given an algebraic group $G$, a maximal torus $T$ of $G$,
 and a Borel subgroup $B$ of $G$ which contains $T$,
the group $\Gamma$ of {\em graph automorphisms} associated to $T$ and $B$,
is defined to be the group of algebraic automorphisms of $G$
 which normalize both $T$ and $B$.  

\begin{fact}[{\cite[\qTheorem 8.4]{BN}}]\label{autalg}
Let $G \rtimes H$ be a group of finite Morley rank,
 with $G$ and $H$ definable.
Suppose that $G$ is a quasisimple algebraic group over an
 algebraically closed field, and $C_H(G)$ is trivial.
Let $T$ be a maximal torus of $G$ and
 let $B$ be a Borel subgroup of $G$ which contains $T$.
Then, viewing $H$ as a subgroup of $\Aut(G)$, we have $H \leq \Inn(G)\Gamma$,
where $\Inn(G)$ is the group of inner automorphisms of $G$ and $\Gamma$
is the group of graph automorphisms of $G$ associated to $T$ and $B$.
\end{fact} 

In other words, quasi-simple algebraic sections never have connected groups
 of outer automorphisms.

\begin{lemma}
Let $H$ be a connected $L$-group of finite Morley rank and odd type
with $\hat O(H)=1$. Let $E_\alg(H)$ be the product of the algebraic
components of $E(H)$.
Then  $$H=E_\alg(H)*H_0$$
where $H_0=C_H\oo(E_\alg(H))$, and $E_\alg(H_0)=1$.
\end{lemma}

\begin{proof}
If $L$ is a component of $E(H)$ with $L/Z(L)$ algebraic,
 then $L$ is algebraic by \cite{AC99}.
Furthermore, a definable group of automorphisms of a quasisimple
algebraic group, in a finite Morley rank context, must itself consist
of algebraic automorphisms, that is graph automorphisms combined with
inner automorphisms, by Fact \ref{autalg}.
So $H$, being connected, induces only inner automorphisms on
$E_\alg(H)$, and this yields
$$ H=E_\alg(H)*H_0 $$
with $H_0=C_H\oo(E_\alg(H))$. 
It is clear that $E_\alg(H_0)\le E_\alg(H)\intersect H_0\le Z(H_0)$,
so $E_\alg(H_0)=1$.
\end{proof}

\begin{theorem}\label{thm:LgrpO}
Let $H$ be a connected $L$-group of finite Morley rank and odd type,
with no normal definable connected solvable subgroup and no normal
definable connected subgroup of degenerate type.
Then $H = H_1 * \cdots * H_n$ is isomorphic to
a central product of quasisimple algebraic groups $H_1,\ldots,H_n$
over algebraically closed fields of characteristic not equal to $2$.
\end{theorem}

\begin{proof}
In this context we have $E(H)=E_\alg(H)$ and $H$ has no nontrivial
definable connected solvable subgroup, so
$F^*(H)=F(H)E_\alg(H)$ with $F\oo(H)=1$. Thus
$C_H(E_\alg(H))=C_H(F^*(H))=Z(F^*(H))$ is finite and $C_H\oo(E_\alg(H))=1$.
So by the preceding lemma, $H=E_\alg(H)$.
\end{proof}

\begin{lemma}\label{Sylow_con}
The Sylow 2-subgroup of a connected $D$-group $G$ of finite Morley rank
is connected.
\end{lemma}

\begin{proof}
We proceed by induction on the rank of $G$.
Suppose first that $G$ has a definable connected normal subgroup $H$
with $1<H<G$.
Let $S$ be a Sylow 2-subgroup of $H$ and
 $S_1$ a Sylow 2-subgroup of $G$ containing $S$.
By \cite[Corollary 1.5.5]{Wag},
 the image of $S_1$ in $G/H$ is a Sylow 2-subgroup of $G/H$,
 and hence connected by induction.
Hence this is also the image of $S_1\oo$ in $G/H$.
So $S_1\le S_1\oo H$, and $S_1=S_1\oo (S_1\intersect H)$.
But $S_1\intersect H$ is a Sylow $2$-subgroup of $H$,
 and is therefore connected.
So $S_1$ is connected.

Now suppose $G$ has no definable connected proper and nontrivial
normal subgroup. If $G$ is abelian, our claim is straightforward
(and more generally if $G$ is solvable it falls under known results).
So we may suppose $Z(G)$ is finite, and then $G/Z(G)$ is simple. 
By the definition of $D$-group, $G/Z(G)$ has degenerate type, and
hence so does $G$. But groups of degenerate type contain no
involutions, so again the Sylow $2$-subgroup is connected.
\end{proof}


\section{$p$-Generation}
\label{sec:Generation}

The proof of our generation result for simple $K^*$-groups and $L^*$-groups 
will be inductive, and depends
therefore on generation results for $K$-groups and $L$-groups. The
latter can be put in a stronger form. The full result runs as follows.

\begin{theorem}\label{generationA}
Let $G$ be a connected group of finite Morley rank, $p$ a prime.
Suppose that $G$ has $p^\perp$ type and $V$ is an elementary abelian $p$-group
 of rank 2 acting definably on  $G$.
Suppose further that one of the following holds.
\begin{enumerate}
\item $G$ is a $K$-group;
\item $G$ is an $L$-group, and $p=2$;
\end{enumerate}
Then $\Gamma_V=G$. 
\end{theorem}

This generation theorem can be strengthened further, but the argument
for this is purely formal. 
Consider the following property.

\begin{definition}
A group $G$ of finite Morley rank has the 
{\em $p$-generation property} if for every elementary abelian
$p$-group $E$ of rank 2 acting definably on $G$, 
we have  $$G=\gen{C_G\oo(v):v\in E^\#}$$
\end{definition}

We observe that this property is trivial when the action of $V$ is not
faithful, so we concern ourselves only with faithful actions.

\begin{lemma}
Let $G$ be a connected group of finite Morley rank such that
every connected definable subgroup of $G$ has the $p$-generation property,
and let $E$ be a nontrivial elementary abelian $p$-group
 of definable automorphisms of $G$.
Then $$G = \gen{C_G\oo(E_0):E_0\le E, [E:E_0]=p}$$
\end{lemma}

\begin{proof}
If $E$ is cyclic this is vacuous. For $E$ of $p$-rank at least 2,
fix $V \le E$ of $p$-rank 2. Then
$$ G = \gen{C_G\oo(v):v\in V^\#} $$
For $v\in V^\#$ let $E = \gen{v} \oplus E_v$. Inductively,
$$ C_G\oo(v) = \gen{ C_G\oo(v,E_1) : E_1\le E_v, [E_v:E_1]=p } $$
and our claim follows.
\end{proof}

Note that in a counterexample to Theorem \ref{generationA}
of minimal Morley rank,
all proper definable connected subgroups and all quotients by 
nontrivial definable connected normal subgroups
will have the $p$-generation property. 

The aim of the present section is to reduce the treatment of Theorem
\ref{generationA} to the simple algebraic case, and to show that the
elementary abelian $p$-groups $V$ involved can be taken to be
contained in the original group $G$.

\subsection{Reductions}

We first treat a ``base case'' of Theorem \ref{generationA}.

\begin{lemma}\label{Gamma-min}
Let $G$ be a connected group of finite Morley rank of $p^\perp$ type
and $V$ an elementary abelian $p$-group of rank 2 consisting of
definable automorphisms of $G$. 
Suppose that one of the following holds.
\begin{enumerate}
\item $G$ is abelian and contains no nontrivial proper definable
  connected $V$-invariant subgroup;
\item  $G$ has degenerate type and $p=2$.
\end{enumerate}
Then $G=\Gamma_V(G)$.
\end{lemma}

\begin{proof}
If $G$ has degenerate type and $p=2$,
 this is \cite[Theorem 5 / Proposition 9.1]{BBC}.

So suppose $G$ is abelian. If $G$ contains no $p$-torsion,
the result is given as  \cite[Fact 3.7]{Bu03}.

Now suppose that $G$ contains $p$-torsion.
Being of $p^\perp$ type, connected, and minimal $V$-invariant,
$G$ is then the definable closure of a $p$-torus $T$.
It suffices to show that $T \leq \Gamma_V$. 

Now $V \le \End(T) \iso M_n(\Zz_p)$ acts on a vector space $W$ over $\Qq_p$
by Fact \ref{endptor}, and the action of $V$ is completely reducible
 in this representation.
Consider an irreducible summand $W_0$ of $W$.
The image of $V$ in $\End(W_0)$ is cyclic, by Schur's lemma.
Thus $W_0$ is centralized by some nontrivial $v \in V$.
In matrix terms, this means that the identity matrix $I$ can be expressed
as a sum of projection matrices $P_i$ with coefficients in $\Qq_p$,
each annihilated by an element corresponding to
the endomorphism $1-v$ for some $v\in V^\#$. 
Thus some multiple $p^n I$ can be similarly expressed in $M_n(\Zz_p)$. 
In $T$ this says that every element of $T = p^n T$ can be expressed
as a linear combination of elements annihilated by $1-v$
for various $v\in V^\#$ and thus $T\le \Gamma_V$. 
\end{proof}

The next result prepares for an inductive argument relative to normal
subgroups. 

\begin{lemma}\label{nonalg:C_v_mod_A}
Let $G$ be a group of finite Morley rank, $p$ a prime,
$v$ a definable automorphism of $G$ of order a power of $p$, and
$K$ a definable normal $p$-invariant connected 
subgroup of $p^\perp$ type. Suppose one of the following holds.
\begin{enumerate}
\item $K$ is solvable.
\item $K$ is of degenerate type and $p=2$.
\end{enumerate}
Then $C_G\oo(v \bmod K) = C_G\oo(v) K$.
\end{lemma}

Here we use the notation $C_G\oo(v \bmod K)$ for the connected
component of the preimage in $G$
of $C_{G/K}(v)$.

\begin{proof}
We may suppose that $G=C_G\oo(v \bmod K)$. Then $G$ is
connected, and $v$ acts trivially on $G/K$.
Furthermore it suffices to treat the case 
in which $v$ has order $p$,
since we can then argue inductively, 
replacing $G$ by $C_G\oo(v^p)$
and $v$ by its restriction to $C_G\oo(v^p)$.

Now in the case $p=2$ with $K$ of degenerate type, 
this is essentially \cite[Lemma 9.3]{BBC} or even \cite[p.~73, Ex.~14]{BN},
together with the fact that connected groups
of degenerate type have no involutions 
(\cite[Theorem 1]{BBC}).

So we may suppose $K$ is solvable.
Proceeding by induction on $\rk(K)$, we may suppose that $K$ contains
no proper nontrivial definable and definably characteristic connected
subgroup.  Then $K$ is abelian.

If $K$ contains no $p$-torsion, then $\gen{v}$ is
 a Sylow $p$-subgroup of $K\gen{v}$, and
 as $[G,v]\includedin K$ we find $K\gen{v}\normal G\gen{v}$.
As Sylow $p$-subgroups are conjugate in $K \gen{v}$,
 the Frattini argument gives $G \leq K\gen{v}\cdot N_G(\gen{v})$. 
But $v$ normalizes $G$, so $N_G(\gen{v}) = C_G(v)$.
Thus $G \le K C_G(v)\gen{v} = K C_{G\gen{v}}(v)$, 
 and therefore $G = K C_G(v)$.

Finally, suppose that $K$ contains nontrivial $p$-torsion. 
By the minimality of $K$,
 as $K$ contains no unipotent $p$-subgroup,
it is the definable closure of a $p$-torus.
As $G$ is connected, we have $K\le Z(G)$.
So the commutation map $h \mapsto [v,h]$ is a
homomorphism from $G$ into $K$, and so has a connected image $L$.
If $v$ centralizes the image $L$ then $L$ has exponent $p$,
 and hence is trivial by connectedness.
Otherwise, our hypothesis on $K$ implies that $C_K(v)$ is finite,
 and thus commutation with $v$ is surjective from $K$ to $K$.
Then for $h\in G$ there is $k\in K$ with $[v,h]=[v,k]$, and
 hence $h k^{-1}\in C(v)$, $h\in KC_G(v)$.
\end{proof}

Now we may carry out an inductive argument.

\begin{lemma}\label{GammaAD}
Let $G$ be a connected group of finite Morley rank of $p^\perp$ type,
for some prime $p$.
Suppose that $G$ contains a nontrivial definable
connected normal definably characteristic subgroup $K$ 
satisfying one of the following conditions.
\begin{enumerate}
\item $K$ is solvable.
\item $K$ has degenerate type, and $p=2$.
\end{enumerate}
If $G/K$ has the $p$-generation property, then so does $G$.
\end{lemma}

\begin{proof}
Fix $V$ a $p$-group of rank 2 consisting of definable automorphisms of
$G$. 
By Lemma \ref{Gamma-min} we have $K\le \Gamma_V$.
By hypothesis
$G=\gen{C_G\oo(v \mod K):v\in V^\#}$, and by the previous
lemma $C_G\oo(v)\le K\Gamma_V(G)$ for $v\in V^\#$.
So putting this all together, $G = \Gamma_V$. 
\end{proof}

We need a variation of this lemma to handle finite normal subgroups.

\begin{lemma}\label{G/Z(G)}
Suppose that $G$ is a connected group of finite Morley rank with
$Z(G)$ finite, $p$ a prime, and $V$ an elementary abelian $p$-group of
rank 2 acting definably on $G$.
If $\Gamma_V(G/Z(G))=G/Z(G)$, 
then $\Gamma_V(G)=G$.
\end{lemma}
\begin{proof}
For $v\in V^\#$, let $G_v=C_G\oo(v\mod\ Z(G))$ 
Then $G=\gen{G_v Z(G)}$ and it suffices to show that $G_v \le \Gamma_V$. 

Now $[v,G_v]\le Z(G)$ is connected, hence trivial. So $G_v \le C\oo(v)\le \Gamma_V$.
\end{proof}

Summing up the reductions so far, we have the following.

\begin{lemma}\label{reduction}
Let $G$ be a connected group of finite Morley rank of $p^\perp$ type
for some prime $p$.
Suppose that every proper definable connected section of $G$ 
has the $p$-generation property, and $G$ does not. 
Then
$F^*(G)=E(G)$ is a direct product of simple groups,
and if $G\ne E(G)$ then $E(G)$ contains no 
$p$-torsion, and $p>2$.
\end{lemma}
\begin{proof}
Fix an elementary abelian $p$-group $V$ of rank 2 acting definably 
on $G$ for which $\Gamma_V < G$.

By Lemmas \ref{GammaAD} and \ref{G/Z(G)} we have
$$Z(G)=1$$
Then by Lemma \ref{GammaAD} again,
we have $F(G)=1$ and $F^*(G)=E(G)$ is a direct product of
connected simple groups. 

Now suppose 
$$E(G)<G$$
Then by hypothesis $E(G)\le \Gamma_V$. 

If $E(G)$ contains $p$-torsion,
 then as $E(G)$ is connected of $p^\perp$ type,
 $E(G)$ contains a nontrivial $p$-torus.
Let $S$ be a Sylow $p$-subgroup of $G V$.
As $G V$ has $p^\perp$ type its Sylow $p$-subgroups are conjugate.
Therefore $S$ contains some maximal $p$-torus $T$ of $E(G)$.
By the Frattini argument $G=E(G)N(T)$, and as $G$ is connected
 we find $$G=E(G)N\oo(T)$$
As $N\oo(T)$ is $V$-invariant and $N\oo(T)<G$
 we also have $N\oo(T) \le \Gamma_V$, so $G = \Gamma_V$, a contradiction. 
So $E(G)$ contains no $p$-torsion.

Now if $p=2$ then $E(G)$ has degenerate type
and we violate Lemma \ref{GammaAD}.
\end{proof}

This amounts to a substantial reduction of the proof of Theorem
\ref{generationA}.
Indeed, if $G$ is a counterexample 
of minimal Morley rank, then $Z(G)$ is finite by
Lemma \ref{GammaAD}, and $G/Z(G)$ is again a counterexample to the
same theorem by Lemma \ref{G/Z(G)}. Then replacing $G$ by $\bar G=G/Z(G)$
we find $Z(\bar G)=1$ and every proper connected section of $G$ has
lower rank, so the previous lemma applies to $\bar G$.

Furthermore, under the hypotheses of Theorem \ref{generationA},
if $E(G)<G$ then either $E(G)$ is a $K$-group, 
hence a product of
simple algebraic groups, and therefore contains $p$-torsion for all
primes $p$, or else $p=2$. In either case Lemma \ref{reduction} is
violated. So in this context, we come down
to the case $G=E(G)$ with $Z(G)=1$.

\subsection{Replacement}

An important tool in the proof of Theorem \ref{generationA}
will be our
ability to replace one elementary abelian $p$-group $V$ by another, as
follows. 

\begin{lemma}\label{replacement}
Let $G$ be a connected centerless group of finite Morley rank 
and $p^\perp$ type for some prime $p$.
Suppose that every proper definable connected subgroup of $G$
 has the $p$-generation property, and
let $E,V$ be two elementary abelian $p$-groups
 of $p$-rank at least two,
acting definably and faithfully on $G$ in
 such a way that the actions commute.
Then $$ \Gamma_E(G)=\Gamma_V(G) $$
\end{lemma}

\begin{proof}
For $v\in V^\#$ we have $C_G\oo(v)<G$ and $E$ acts on $C_G\oo(v)$, so
by hypothesis $C_G\oo(v)\le \Gamma_E$, and hence $\Gamma_V \le \Gamma_E$.
Similarly $\Gamma_E\le \Gamma_V$.
\end{proof}

We note that if $G$ is both of $p^\perp$-type and $q^\perp$ type for
two primes $p,q$ then the same result holds with $V$ 
an elementary abelian $p$-group and $W$ an elementary abelian
$q$-group. This can be useful with $p>2$ and $q=2$.

It is natural at this point to look at the graph $\EE_{p,2}$ whose
 vertices are the elementary abelian $p$-groups acting definably on $G$
(and for definiteness, living in the normalizer of $G$ in
some larger group in which the centralizer of $G$ is trivial), 
and whose
edges are commuting pairs of such groups. Then $\Gamma_V$ is
associated to a connected component of this graph, rather than to a
single vertex. This graph is frequently connected, in which case
$\Gamma_V$ is a canonical normal subgroup of $G$---which is
not surprising, as we expect it 
to be the group $G$ itself in the cases of interest to us.
We will write $\EE_{p,2}(X)$ for the corresponding graph with vertices
contained in $X$ where $X$ is a fixed group normalizing $G$.

The following is \cite[Fact 1.20]{Bu07a}.

\begin{fact}\label{p+1}\label{Ep2}
Let $S$ be a locally finite $p$-group.  Then
\begin{conclusions}
\item The induced subgraph $\EE^0_{p,2}(S) :=
  \{ X \in \EE_{p,2}(S) : m_p(C_S(X)) > 2 \}$ is connected, i.e., 
  all nonsingleton vertices of $\EE_{p,2}(S)$ are linked by paths.
\item $\Omega_1(S\oo)$ is found in $\EE^0_{p,2}(S)$.
\item If $S$ contains an abelian normal subgroup of $p$-rank at least $p+1$
  then $\EE_{p,2}(S)$ is connected.
\end{conclusions}
\end{fact}

This has the following immediate consequence.

\begin{lemma}
Let $G$ be a connected group of finite Morley rank of $p^\perp$ type
 for some prime $p$, which is a direct product of connected simple groups.
Suppose that $G$ does not have the $p$-generation property, but that
 every proper definable connected subgroup of $G$ does.
Then one of the following holds:
\begin{enumerate}
\item $G$ is simple;
\item $G$ is the product of exactly $p$ isomorphic simple components, 
and each component $L$ has an automorphism of order $p$ which does not
normalize any nontrivial proper definable connected subgroup of $L$. 
In particular, $G$ contains no $p$-torsion.
\end{enumerate}
\end{lemma}

\begin{proof}
Let $V=\gen{a,b}$ be an elementary abelian $p$-group
 of $p$-rank 2 acting definably on $G$ so that $\Gamma_V<G$. 
Suppose $G$ contains a proper nontrivial normal $V$-invariant
subgroup $H$. Then $G=H \times C_G(H)$ and by hypothesis both $H$ and
$C_G(H)$ are contained in $\Gamma_V$, so $G=\Gamma_V$, a
contradiction.

Accordingly $V$ acts transitively on the simple components of $G$. 
Suppose first that $V$ acts regularly on the set of components.
Let $L$ be a component of $G$, and let $H,K$ be the
normal closures of $L$ under the action of $\gen{a}$ and $\gen{b}$
respectively. Then $C_H(a)$, $C_K(b)$ are diagonally embedded copies
of $L$ in $H$ and $K$ respectively, and $H\intersect K=L$. 
It follows that $[C_H(a),C_K(b)]=L$. Thus $L\le \Gamma_V$ and so
$G=\Gamma_V$, a contradiction.

So $V$ does not act regularly on the components of $G$ and hence we
may suppose the component $L$ of $G$ is normalized by $a$. 
If $L$ is normalized by $V$ then $G=L$ is simple. So we remain with
the case in which $G$ is the product of precisely $p$ conjugates of
$L$. Now suppose $a$ normalizes a nontrivial proper definable 
connected subgroup $L_0$ of $L$. Then the product of the conjugates of
$L_0$ under $V$ lies in $\Gamma_V$ by hypothesis, and in particular
$L_0\le \Gamma_V$.
On the other hand,
$C_G\oo(b)$ is a diagonal subgroup of $G$ covering $L$ in this case
and $[C_G\oo(b), L_0]=L$ and so $L\le \Gamma_V$ and then $G=\Gamma_V$, a
contradiction. 

So $a$ normalizes no nontrivial proper definable connected subgroup of
$L$. In particular $L$ contains no $p$-torsion, as otherwise $L$ would
contain an $a$-invariant maximal $p$-torus, using the conjugacy of 
Sylow $p$-subgroups in $L\gen{a}$.
So $G$ contains no $p$-torsion. 
\end{proof}

\begin{corollary}
For the proof of Theorem \ref{generationA}, it suffices to treat the
case of simple groups $G$, under the hypothesis that every proper
definable connected section of $G$ has the $p$-generation property.
\end{corollary}

\begin{proof}
Indeed, we have already seen that after minimization and factoring out the
center, the group $G$ will be a direct product of connected simple groups.
If $G$ is not simple, then
 the previous lemma shows that $G$ contains no $p$-torsion.
On the other hand, under the hypotheses of Theorem \ref{generationA},
 $G$ is either a $K$-group, or an $L$-group with $p=2$.
In the former case as the factors of $G$ are algebraic,
$G$ contains $p$-torsion for all $p$, and we contradict the lemma.
In the latter case as $p=2$ and $G$ contains no $p$-torsion,
 $G$ is of degenerate type and Lemma \ref{Gamma-min} applies.
\end{proof}

Finally, we argue that one may move the elementary abelian $p$-groups
under consideration inside the group $G$, when $G$ is algebraic.

\begin{lemma}\label{inner_reduction}
Let $G$ be a simple algebraic group of finite Morley rank of $p^\perp$ type
such that all proper definable connected sections of $G$ have
 the $p$-generation property.
Let $V$ be an elementary abelian $p$-group of rank two
 acting definably and faithfully on $G$, with $\Gamma_V<G$.
Then there is an elementary abelian $p$-group $U$ of rank two
 contained in $G$ such that $\Gamma_V=\Gamma_U$.
\end{lemma}

\begin{proof}
By Fact \ref{p+1}, if the Lie rank of $G$ is at least 
$p+1$ and $S$ is a Sylow $p$-subgroup of $G V$, then the graph 
$\EE_{p,2}(S)$ is connected. Then by Lemma \ref{replacement} 
we have $\Gamma_V=\Gamma_U$ for any $U\le \Omega_1(S\oo)$ of rank two.

On the other hand, if $G$ has no graph automorphisms of order $p$ then
any definable automorphism of $G$ of order $p$ will be inner,
 by Fact \ref{autalg}, in which case we have nothing to prove.
As the graph automorphisms have orders $2,3$ we consider the two cases $p=2,3$.

If $p=3$ we are dealing with the group $D_4$ which has Lie rank $p+1$,
and as noted we are done in this case.

If $p=2$ we need to deal only with groups of Lie rank at most two 
and the only one with a graph automorphism is $A_2$, namely $\PSL_3$. 
In this case a graph automorphism $v\in V$ is conjugate to the
 inverse-transpose automorphism, and has centralizer $C_G(v) \cong \PSL_2$
 by \cite[Table 4.3.1 p.\ 145 \& Table 4.3.3 p.\ 151]{GLS3}.
On the other hand $V$ must also contain an element $u$ inducing
an inner automorphism; taking $u$ without loss of generality to be an
element of $G$ we have $u\in C_G(v)$. Now $u$ belongs to a 4-group
$U$ in $C_G(v)$, and then  $U$ commutes with $V$, so by Lemma
\ref{replacement} $\Gamma_V=\Gamma_U$.
\end{proof}

\section{$p$-Generation in simple algebraic groups}\label{sec:SAG}

Our goal now is to prove Theorem \ref{generationA} for simple 
algebraic groups over an algebraically closed field, using the
machinery of linear algebraic groups \cite{Springer_LinAlgGrp,Hump},
as well as more detailed data on simple algebraic groups \cite{Carter93}.
However, the reader may recognize that some results used here are
merely facts about groups of finite Morley rank applied in the reduct
to the pure field language.

As we have seen, we may suppose that all proper definable connected
subgroups satisfy $p$-generation, and we may restrict our attention
to $\Gamma_V$ for $V$ contained in the given group by Lemma \ref{inner_reduction}.
So we normally work under the following hypotheses, with $p^\perp$ type
translated to characteristic not $p$.

\begin{hypothesis}\label{hyp:pgeneration_inductive}
Let $G$ be an algebraic group over an algebraically closed field
 of characteristic other than $p$.
Suppose that every proper connected
 Zariski closed subgroup has the $p$-generation property.
\end{hypothesis}

It will however be convenient to dispose of $\PSL_2$ at the outset.

\begin{lemma}
If $G = \PSL_2$ over an algebraically closed field of characteristic
other than $p$, then $G$ satisfies $p$-generation.
\end{lemma}

\begin{proof}
We are dealing with an elementary abelian $p$-group $V$ of rank 2
with $V\le G$. Each element of $V$ belongs to a maximal torus of $G$ 
 by \cite[6.3.5]{Springer_LinAlgGrp} or Fact \ref{torality}.
If these tori do not lie in a Borel subgroup they generate $G$.
So $\Gamma_V = G$.
But if they do lie in a Borel subgroup then so does $V$,
contradicting the structure of Borel subgroups in this case.
\end{proof}

Now we undertake some structural analysis.  The following extreme case
arises frequently.

\begin{lemma}\label{Borel}
Suppose that $G$ satisfied Hypothesis \ref{hyp:pgeneration_inductive}.
Let $V$ be an elementary abelian $p$-subgroup of rank two contained
 in a Borel subgroup of $G$.
Then $\Gamma_V = G$.
\end{lemma}

\begin{proof}
As $G$ is not $\PSL_2$, the maximal parabolic subgroups of $G$
containing $B$ generate $G$, and by hypothesis each such subgroup
is contained in $\Gamma_V$.  So $G=\Gamma_V$.
\end{proof}

\begin{lemma}\label{reductive}
Let $G$ be a simple group satisfying Hypothesis \ref{hyp:pgeneration_inductive}.
Then $\Gamma_V$ is a Zariski closed reductive subgroup of $G$
 containing a maximal torus of $G$, and with finite center;
 the center has order relatively prime to $p$.
\end{lemma}

\begin{proof}
Every element of $V$ lies in a maximal torus of $G$ and thus 
$\Gamma_V$ contains a maximal torus of $G$.
Let $H$ be the subgroup of $\Gamma_V$ generated by such maximal tori of $G$,
 and suppose that $H < G$.
Then $H$ is Zariski closed and $\Gamma_V\le N\oo(H)<G$.
By hypothesis, $\Gamma_V=N\oo(H)$, a Zariski closed subgroup.

If $U$ is the unipotent radical of $\Gamma_V$ then for $v\in V^\#$ 
we have $U\intersect C\oo(v)$ contained in the unipotent radical of
$C\oo(v)$; but $v$ being a semisimple element, $C\oo(v)$ is reductive
\cite[Thm.~3.5.4]{Carter93}, and thus $C_U\oo(v)=1$. But by hypothesis
$U=\gen{C_U\oo(v):v\in V^\#}$ and thus $U=1$, $\Gamma_V$ is reductive.

Now we claim that $Z(\Gamma_V)$ contains no element of order $p$.
As $G$ has characteristic not $p$, it then follows that $Z(\Gamma_V)$
contains no torus, and hence $Z(\Gamma_V)$ is finite of order prime to $p$.
So this will complete the analysis.

So suppose toward a contradiction that 
$u\in Z(\Gamma_V)$ has order $p$. 
We can replace $V$ by a subgroup of $\gen{V,u}$
containing $u$. So we may suppose $u\in V$, say $V=\gen{u,v}$.
Now $v$ belongs to a $p$-torus $T$ of $\Gamma_V$ and $u$ commutes with $T$. 
Take a maximal $p$-torus $T_1$ of $C(u)$ containing $T$.
Then $v\in T_1$ by definition and $u\in T_1$ by Fact \ref{torality}.  
Thus $V$ is contained in a torus of $G$, and in particular is contained
in a Borel subgroup of $G$.  So $\Gamma_V = G$ by Lemma \ref{Borel},
 and $Z(\Gamma_V)=1$, a contradiction.
\end{proof}

Now Lemma \ref{reductive} puts us in a position to control $\Gamma_V$
even more tightly.

\begin{lemma}\label{quasisimple}
Let $G$ be a simple group satisfying Hypothesis \ref{hyp:pgeneration_inductive},
 and let $V\le G$ an elementary abelian $p$-subgroup of $p$-rank two. 
Then $\Gamma_V$ is quasisimple.
\end{lemma}

\begin{proof}
We know that $\Gamma_V$ is a central product of quasisimple factors,
and that the center of $\Gamma_V$ contains no elements of order
$p$. Therefore the Sylow $p$-subgroup of $\Gamma_V$ is the direct
product of the Sylow $p$-subgroups of the factors.

Now as every element of $V$ belongs to a torus by Fact \ref{torality},
 we have $V\le \Gamma_V$ and
 thus $V$ is contained in a Sylow $p$-subgroup of $\Gamma_V$. 
Let $V=\gen{u,v}$ and choose quasisimple components $L_u,L_v$ of $G$
so that the projection $u'$ of $u$ into $L_u$ and $v'$ of $v$ into $L_v$
is nontrivial. Let $U=\gen{u',v'}$.
Then $U,V$ commute and we may replace $V$ by $U$.
Evidently $\Gamma_U=G$ via Lemma \ref{Borel} unless $L_u=L_v$.
In this case writing $L=L_u=L_v$, we have $U\le L$ and
if $L < \Gamma_V$ then we may replace $U$ by a subgroup
$<u',v''>$ with $v''$ coming from another component of $\Gamma_V$.
We are left with the possibility $L=\Gamma_V$.
\end{proof}

Now we require some understanding of
 the universal cover of a simple algebraic group.

\begin{lemma}\label{algCindex}
Let $\hat{G}$ be the universal central extension of
a quasisimple algebraic group $G$.
Then for any semisimple element $s\in G$ we have
$$\hbox{$[C_G(s) : C_G\oo(s)]$ divides  $|Z(\hat{G})|/|Z(G)|$}$$
\end{lemma}

\begin{proof}
A maximal torus of $\hat G$ covers a maximal torus of $G$.
So there is a semisimple element $t\in \hat G$ corresponding to $s$
 under the quotient map $\hat G\to G$.
In simply connected quasisimple algebraic groups the centralizers of 
semisimple elements are connected (\cite[Theorem 3.5.6]{Carter93}).
Thus $C_{\hat G}(t)$ is connected.

Of course $Z := Z(\hat G)$ is the kernel of this map $\hat G\to G$. 
Now $K=C_G(s)$ pulls back to the subgroup $\hat K = C_{\hat G \mod\ Z}(t)$
of $\hat G$.  Taking the commutator with $t$ we get a homomorphism
$\hat K \to Z$ with kernel $C_{\hat G}(t)$.
As $C_{\hat G}(t)$ is connected, we have $\hat K\oo = C_{\hat G}(t)$
and $|\hat K/\hat K\oo|$ divides $|Z|$.

Now as $Z(\hat G)\le \hat K\oo$ we have $|K/K\oo|=|\hat K/\hat K\oo|$
and our claim follows.
\end{proof}

The Schur multiplier of a simple algerbaic group is the center of its
(algebraic) universal central extension.

\begin{lemma}\label{Schur}
Let $G$ be a group satisfying Hypothesis \ref{hyp:pgeneration_inductive}.
Suppose that $p$ does not divide the order of the Schur multiplier.
Then $G$ has the $p$-generation property.
\end{lemma}

\begin{proof}
As we have shown this comes down to the case of an elementary abelian
$p$-group of $p$-rank two contained in $G$.
We will show that $V$ is contained in a torus of $G$.
Write $V=\gen{u,v}$ and consider $C\oo(u)$. 
By Lemma \ref{algCindex}, the index $[C(u):C\oo(u)]$ is not divisible
by $p$ and therefore $v\in C\oo(u)$.
Now $v$ belongs to a maximal torus $T$ of $C\oo(u)$ and
 $u$ also belongs to $T$ by \cite[6.3.5]{Springer_LinAlgGrp} or Fact \ref{torality},
so $V$ is contained in a torus of $G$.
\end{proof}

The orders of the Schur multipliers of simple algebraic groups are
 given in \cite[\S1.11, p.\ 25--26]{Carter93}.
Namely: $n+1$ for $A_n$ with $n\geq1$;  either 2 or 4 for
 each of $B_n$, $C_n$, $D_n$ with $n\geq2,2,3$ respectively, and $E_7$;
 3 for $E_6$; and finally, $1$ for $G_2$, $F_4$, and $E_8$.

\begin{proof}[Proof of Theorem \ref{generationA}]
After considering a minimal case, factoring out the center, and using
our replacement lemma, we arrive at the following situation:
\begin{enumerate}
\item $G$ is a simple algebraic group of $p^\perp$ type.
\item $V\le G$ is an elementary abelian $p$-subgroup of $p$-rank two.
\item Every proper definable connected section of $G$ has the
  $p$-generation property.
\item $\Gamma_V<G$.
\end{enumerate}

Furthermore $\Gamma_V$ contains a maximal torus of $G$ and is
quasisimple. By considerations of dimension this already eliminates
the possibility that $G$ has type $A_n$ ($\PSL_{n+1}$).

Furthermore the prime $p$ involved has to divide the order of the
Schur multiplier of $G$, and by the calculations of Schur multipliers
in the remaining cases this brings us down to the cases $p=2,3$.

On the other hand, by Fact \ref{p+1}, if the Lie rank of $G$ is at
least $p+1$ we can replace $V$ by a subgroup of a torus and conclude
by Lemma \ref{Borel}. So looking at the orders of the Schur
multipliers again 
we eliminate the case $p=3$ and for $p=2$ we come down to type $B_2$.
But then as $\Gamma_V$  is itself quasisimple of Lie rank two,
the only possibility compatible with the dimensions would be an
embedding of $\SL_3$ into $B_2$, which does not occur because the orders of
their Weyl groups are given by $|A_n| = (n+1)!$ and $|B_n| = 2^n n!$ with $n=2$.
\end{proof}
\section{Generation in $K^*$-groups and $L^*$-groups}

We turn now to the proof of our main generation results.

\begin{lemma}
Let $G$ be a connected simple group of finite Morley rank of $p^\perp$ type.  
Suppose that any proper definable connected subgroup of $G$
 has the $p$-generation property and that
 $V$ is an elementary abelian $p$-subgroup of $G$ of rank two.
Let $S$ be a Sylow $p$-subgroup of $G$ containing $V$,
 and suppose that $m_p(C_S(V))\ge 3$. 
Then
\begin{conclusions}
\item $N(\Gamma_V)=\Gamma_{S,2}^0$
\item $\Gamma_V=(\Gamma_{S,2}^0)\oo$.
\end{conclusions}
\end{lemma}

\begin{proof}
Let $M=N(\Gamma_V)$. Then $V$ normalizes $M\oo$.
If $M\oo<G$ then $M\oo$ satisfies $p$-generation and hence $M\oo\le
\Gamma_V$. If $M\oo=G$ then $\Gamma_V\normal G$ and hence
$\Gamma_V=G$. So in either case $M\oo\le \Gamma_V$.
Thus 
$$M\oo=\Gamma_V$$

We show next that $$ \Gamma_V \le \Gamma_{S,2}^0 $$
We take $v\in V^\#$, take $E\le S$ to be an elementary abelian $p$-group
of rank three containing $V$, and write $E = \gen{v}\oplus E_v$ with
$E_v$ elementary abelian of rank two containing $v$.
Then by $p$-generation we have
$$C\oo(v)=\gen{C\oo(\gen{v,w}):w\in E_v}$$
and as $C\oo(v,w)\le N(\gen{v,w})\le \Gamma_{S,2}^0$ for such pairs
$v,w$ we find $\Gamma_V\le \Gamma_{S,2}^0$.

On the other hand by Fact \ref{Ep2} we have $\Gamma_V=\Gamma_W$ for
$W\le \Omega_1(S)$ and in particular we find $S\le \Gamma_V$. 
So by the Frattini argument $N(\Gamma_V)\le \Gamma_V\cdot N(S\oo)$ and
as $N(S\oo)\le N(\Omega_1(S\oo))\le \Gamma_{S,2}^0$ we find
$$N(\Gamma_V)\le \Gamma_{S,2}^0$$

Now we show $\Gamma_{S,2}^0\le N(\Gamma_V)$. We take $W\le S$ an
elementary abelian $p$-group with $m_p(C_S(W))\ge 3$ and observe that
$\Gamma_V=\Gamma_W$ by Lemma \ref{replacement}. Therefore $N(W)\le
N(\Gamma_V)$ for all such $W$, and $\Gamma_{S,2}^0\le N(\Gamma_V)$. 
With this all claims are proved.
\end{proof}

\begin{corollary}\label{cor:2-gencore}
Let $G$ be a connected simple group of finite Morley rank of $p^\perp$
type with $m_p(G)\ge 3$.  Suppose that any proper definable
connected subgroup of $G$ has the $p$-generation property.
Then the following are equivalent.
\begin{enumerate}
\item For some elementary abelian $p$-subgroup $V$ of $G$ of $p$-rank two, 
$\Gamma_V<G$.
\item For $S$ a Sylow $p$-subgroup of $G$, we have $\Gamma_{S,2}^0<G$.
\end{enumerate}
\end{corollary}

\begin{proof}
By the conjugacy of Sylow $p$-subgroups
 and Theorem \ref{max_el_abelian},
every elementary abelian $p$-subgroup $V$ of $p$-rank two in $G$
 is contained in a Sylow $p$-subgroup $S$ satisfying $m_p(C_S(V))\ge 3$,
and every Sylow $p$-subgroup $S$ of $G$ contains such
 an elementary abelian $p$-subgroup $V$.
Since $N(\Gamma_V)=\Gamma_{S,2}^0$ under these conditions, our claim follows.
\end{proof}

Now we state the key result.

\begin{theorem}\label{thm:strongembedding}
Let $G$ be a connected simple group of finite Morley rank of $p^\perp$ type 
 with $m_p(G)\ge 3$ and $m_2(G)\ge 2$.
Suppose that any proper definable connected subgroup of $G$ has
 the $p$-generation property and the 2-generation property.
Suppose that $V$ is an elementary abelian $p$-subgroup of $G$ and $\Gamma_V<G$.
Then $N(\Gamma_V)$ is strongly embedded in $G$.
\end{theorem}

By Theorem \ref{generationA}, our condition on the sections of $G$
will hold if $G$ is a $K^*$-group, or if it is an $L^*$-group with
$p=2$, giving the two cases of concrete interest.

We record a standard criterion for strong embedding
(\cite[Thm.~10.20]{BN}). 

\begin{fact}\label{strembdef}
A subgroup $M$ of a group $G$ is strongly embedded if and only if $M$
satisfies the following two conditions.
\begin{enumerate}
\item $M$ contains the normalizer of a Sylow $2$-subgroup of $G$;
\item $M$ contains $C(i)$ for each involution $i\in M$.
\end{enumerate}
\end{fact}

\begin{proof}[Proof of Theorem \ref{thm:strongembedding}]
Let $M=\Gamma_V$.

If $G$ is not of odd type then it is algebraic of even type by \cite{ABC_EvenType}.
In particular, $G$ is a $K$-group in the pure field language.
It follows that $\Gamma_V = G$ by Theorem \ref{generationA}.
So we assume $G$ has odd type.

Now we may suppose $V\le S$ with $S$ a Sylow $p$-subgroup 
chosen so that ${m_p(C_S(V))\ge 3}$.

We show first that $\Gamma_V$ contains a Sylow\oo\ 2-subgroup of $G$. 
Let $T$ be a maximal divisible abelian torsion subgroup of $G$, and
$T_p$, $T_2$ the $p$-torsion and 2-torsion subgroups of $T$.
The maximal divisible abelian torsion subgroups of $G$ are conjugate
 \cite{ABC_EvenType,Ch05} and thus $T_p$ and $T_2$ are respectively
 a maximal $p$-torus and a maximal $2$-torus of $G$.
We may suppose after conjugating that $T_p\le S$.
Take $A\le \Omega_1(T_p)$ and $B\le \Omega_1(T_2)$
 elementary abelian subgroups of rank 2.
Then by Fact \ref{Ep2} and Lemma \ref{replacement} we have
$\Gamma_V=\Gamma_A$. By the variant of Lemma \ref{replacement} 
mentioned after its proof, we also have $\Gamma_A=\Gamma_B$. 

Embed $T_2$ in a Sylow 2-subgroup $P$ of $G$. Then $P\oo=T_2$ and
$N(P) \le N(\Omega_1(T_2)) \le N(\Gamma_B)=N(\Gamma_V)$. Thus $N(P)\le M$.

Next we show that for $i$ an involution of $M$, we have $C\oo(i)\le
M$. We may suppose that $i$ normalizes $T_2$. If $i$ centralizes $B$,
then $i$ lies in an elementary abelian $2$-subgroup $U$ of $2$-rank 
two that commutes
with $B$, and hence $\Gamma_V=\Gamma_B=\Gamma_U$ and $C\oo(i)\le
\Gamma_U\le M$. So suppose that $i$ does not centralize $B$ 
and hence does not
invert $T_2$. Therefore $i$ centralizes a nontrivial $2$-torus $T_0\le
T_2$ by Lemma \ref{centinvT}.
By Fact \ref{torality} if $\hat T_0$ is a maximal $2$-torus in
$C(i)$, then $i\in \hat T_0$. Now $\hat T_0\le \Gamma_V$ since $\hat
T_0$ contains an involution of $T_2$.
Hence after
conjugation in $\Gamma_V$ we may suppose $\hat T_0\le T_2$ and thus $i\in
T_2$. But then as above our claim follows.

Now for $i$ an involution of $M$ it remains to show that $C(i)\le M$.
By Fact \ref{torality} the element $i$ lies in some maximal $2$-torus 
$T_i$ of $G$, and then $T_i\le C\oo(i)\le \Gamma_V$ and after
conjugation we may suppose $T_i=T_2$. Then by the Frattini argument
$$C(i)\le C\oo(i)\cdot N(T_2)$$
and as $N(T_2)\le N(\Gamma_V)=M$ our claim follows.
\end{proof}

\section{Strongly Embedded Subgroups}
\label{sec:Altseimer}

We next treat the structure of strongly embedded subgroups of
$L^*$-groups, by a method used by Altseimer \cite{AltseimerStrEmb} in
the $K^*$ context.

\begin{theorem}\label{thm:Altseimer}
Let $G$ be a simple $L^*$-group of finite Morley rank and odd type.
Then a strongly embedded subgroup $M$ of $G$ is
either a $D$-group, or else satisfies $M^\o/\hat O(M^\o) \cong \SL_2$.
\end{theorem}

We begin by recalling an elementary property of groups with strongly
embedded subgroups. 

\begin{fact}[{\cite[Theorem 10.19]{BN}}]
\label{stremb_invconj}
Let $G$ be a group of finite Morley rank with a proper definable
strongly embedded subgroup $M$.  Then 
$G$ and $M$ each have only one conjugacy class of involutions.
\end{fact}

For the proof, one shows first that any involution in $M$ is conjugate
to any involution not in $M$, after which the conjugacy in $G$ is
clear and the conjugacy in $M$ follows from the definition of strong
embedding. 

\begin{lemma}\label{Lgrp_oneconjclass}
Let $M$ be an $L$-group of finite Morley rank and odd type 
with exactly one conjugacy class of involutions.
Let $H=M^\o/\hat O(M)$.
Then either 
\begin{enumerate}
\item $M$ is a $D$-group, and $H/Z(H)$ 
  contains no involutions; or else 
\item $H$ is quasisimple, and
 and $H/Z(H)$ is one of the 
 following: $\PSL_2$, $\PSL_3$, or $\G_2$. 
\end{enumerate}
\end{lemma}

\begin{proof}
We need a few facts about algebraic groups.
First, it is well known, and can be extracted from Theorem 1.12.5d of \cite{GLS3},
that $\SL_2$ is the only quasisimple algebraic group whose involutions
all lie in its center. 
Second, as can be seen in Table 4.3.1 on p.\ 145 of \cite{GLS3},
any algebraic group which is simple as an abstract group, and
 has a unique conjugacy class of involutions,
must be one of the following: $\PSL_2$, $\PSL_3$, or $\G_2$.
All others contain involutions with nonisomorphic centralizers, and
hence have distinct conjugacy classes of involutions that do not even fuse
under the action of $\Aut(G)$.

Now suppose first that 
$$\hbox{$H$ contains no nontrivial connected normal
definable $D$-subgroup.}\leqno(I)$$ 
Then by Theorem \ref{thm:LgrpO} $H = E_\alg(H)$.

If $Z(H)$ contains an involution, then
 all involutions in $H$ are central, and
 the components are all of the form $\SL_2$.
One may easily verify that
 $\SL_2 \times \SL_2$ has one conjugacy class of non-central involutions, and
 a central product $\SL_2 * \SL_2$ with $Z = 1$ has two.
So it follows that there is a unique component and $H = \SL_2$
 because all involutions lie in $Z(H)$, 
We may now put this case aside since $H \iso \SL_2$ is one of our conclusions.

So suppose $Z(H)$ contains no involution. Then each component contains
noncentral involutions and hence every involution lies in a unique
component. This again implies that there is only one component. 
Furthermore, all of its involutions are conjugate. So $H$ is again one
of those listed.

Now suppose 
$$\hbox{$H$ contains some nontrivial connected normal
definable $D$-subgroup $K$.}\leqno(II)$$ 

We may take $K$ minimal. Note that $\hat O(K)=1$.
If $F^\o(K)=1$ then $E(K)$ must be nontrivial, and for each quasisimple
component $L$ of $E(K)$ we have $L/Z(L)$ of degenerate type. But also
$Z(L)\le F(K)$ is finite, and thus $L$ is of degenerate type, and
connected. So $L$ contains no involutions and $L\le \hat O(K)=1$, a
contradiction. Accordingly $$ F^\o(K) > 1 $$
Again, $O(F^\o(K))=1$, so the 2-torsion in $F^\o(K)$ forms a
 nontrivial 2-torus $T$.
As $N^\o(T)=C^\o(T)$, this 2-torus is central in $H$.
It follows that all 2-elements of $H$ lie in $T$
 by conjugacy of maximal decent tori \cite{Ch05},
Thus $H/Z(H)$ contains no involutions, and $H$ is a $D$-group.
\end{proof}

This moves us substantially in the direction of Theorem \ref{thm:Altseimer}.
But we still need to eliminate the configurations in which
$M^\o/\hat O(M)$ is 
$\PSL_2$, $\PSL_3$, $\SL_3$, or $\G_2$.
For this we must further exploit strong embedding in the manner of
\cite[\S5]{AltseimerStrEmb}, and also (in a different direction) \cite[\S5]{AC03}.

A powerful idea in this context is to study the distribution of
involutions in cosets of $M$ (or $M^\o$). This leads to the following.

\begin{fact}[{\cite[Lemma 3.8]{Al96}}]
Let $G$ be a group of finite Morley rank
 with a strongly embedded subgroup $M$.
Then there is an involution $w\in G \setminus M$
 such that $\rk(I(w M)) \geq \rk(I(M))$.
\end{fact}

Fix such an involution $w$, set $Y := \{ u w \mid u \in I(w M) \}$,
 and let $K := d(Y)$.

\begin{fact}[{see \cite[Propositions 3.9,3.10]{Al96}}]
\label{stremb_K_no_inv}\label{stremb_K_plus_Ci}
\label{CK}
Let $G$ be a group of finite Morley rank
with a strongly embedded subgroup $M$.
Then the group $K$ contains no involutions,
and for any involution $i\in M$, we have 
$$M^\o = C^\o(i) K^\o$$
\end{fact}

This decomposition will then pass to the quotient $M^\o/\hat O(M)$, or for
that matter a further quotient by the center, if the center has odd order.
So we require the following.

\begin{lemma}\label{CKanalysis}
Let $H$ be one of the groups $\PSL_2$, $\PSL_3$, or $\G_2$,
over an algebraically closed field of characteristic not $2$.
Then there is no connected definable subgroup $K$ of $H$ without
involutions for which we have a decomposition
$$ H = C(i) K $$  for each involution $i\in H$.
\end{lemma}

The following general fact will be useful in our analysis.

\begin{lemma}\label{algebraicity_of_unipotence}
Any definable unipotent subgroup $U$ in an algebraic group $H$
over a field $k$ of finite Morley rank and characteristic zero is Zariski closed.
\end{lemma}

We note that the language in which $U$ is supposed to be definable
(the language of $k$) is arbitrary, subject to the finite rank hypothesis.

\begin{proof}
The additive group of a field $k$ of finite Morley rank and characteristic zero
has no proper (infinite) definable subgroup \cite[Cor.~3.3]{Po87}.
It follows that a definable endomorphism $\alpha$ of $(k,+)$ is
given by multiplication in $k$. 

Since any such group $U$ would be nilpotent, 
we may assume, proceeding inductively, that $U$ is abelian, and minimal.
Let $\hat{U}$ be the Zariski closure of $U$.
Then $\hat U$ can be identified with a vector space over $k$
(cf.~\cite[Lemma 15.1C]{Hump}).
Consider a projection into any $1$-dimensional subspace $V$ of $\hat U$.
Since $\hat U$ is the Zariski closure of $U$, $U$ is not contained in
the kernel. Hence $U$ maps injectively into $V$. 
So the projection gives an isomorphism.
Writing $\hat U=\oplus_i V_i$, we can identify the $V_i$ with $k$ and
view $U$ as being given by a series of definable automorphisms of $(k,+)$.
These are represented by multiplication maps on $k$ and thus $U$ is
defined by linear equations.
\end{proof}

\begin{proof}[Proof of Lemma \ref{CKanalysis}]
We make use of rank computations for algebraic subgroups of $H$. If
$r$ is the rank of the base field, the rank of such a subgroup is just
$r$ times its dimension. So we may work directly with dimensions.

Suppose first that $K$ is solvable, and let $B$ be a Borel subgroup of
$H$ containing $K$. (We remark that the model theoretic and algebraic
notions of Borel subgroup coincide here, but in any case we will need
the structure theory of algebraic groups, so the algebraic sense is
dominant.)  Taking $i$ in $B$, and discounting the overlap caused by
a maximal torus $T$ of $B$ containing $i$, we find
$$\dim(H)\le \dim(C(i))+\dim(B)-\dim(T)$$
where one may replace $\dim(B)-\dim(T)$ by $\dim(U)$ with $U$ the
unipotent radical of $B$, so that this term is the number of positive
roots. The relevant numbers are found in the following table.

\begin{table}[h] 
\begin{center}
\begin{tabular}{cccccc} \toprule
$H$ & Lie rank & $\dim(C(i))$ & $\dim(U)$ &Sum & $\dim(H)$ \\ \midrule
$\PSL_2$ & 1 &  1 & 1 &2&3\\
$\pPSL_3$ & 2 & 4 & 3 &7&8\\
$\G_2$ & 2 & 6 & 6 &12&14\\
\bottomrule
\end{tabular}
\end{center}
\caption{Dimensions of some subgroups}\label{T:alg_subgrp_dim}%
\end{table}
\noindent Looking at the last two columns, we see that our conditions are not
met. 
So we conclude
$$\hbox{The group $K$ is not solvable.}$$

By a subtle result of Poizat, if the base field for $H$ has positive
characteristic then every connected definable simple section of $H$ 
is algebraic \cite[Theorem 1]{Po01a} (cf.~\cite[4.19]{AC03} for the form quoted here).
This then forces $K$ to contain involutions and gives a contradiction.
So the base field has characteristic $0$.

Now let $\hat K$ be the Zariski closure of $K$.
Observe that $\hat K$ inherits from $H$ the key property:
$$\hbox{ If $i$ is an involution in $\hat K$ then $\hat K = C_{\hat K}(i)K$. }$$
Furthermore $\hat K/U(\hat K)$ (factoring out the unipotent radical)
inherits this property. This is a reductive group. There are two
possibilities: either it is quasisimple of Lie rank $2$, or every
quasisimple component, modulo its center, is $\PSL_2$.

In the latter case, look at the image of $K$ in $\hat K/U(\hat K)$ and
then projected into the simple components modulo the center. These
give subgroups of $\PSL_2$ containing no involutions, which are
solvable by another result of Poizat's \cite[Theorem 4]{Po01a}.
In this case $K$ itself is solvable and we have a contradiction.
Thus $$\hbox{$\hat K/U(\hat K)$ is quasisimple of Lie rank $2$.}$$

Replacing $H$ by $\hat K/U(\hat K)$ and $K$ by its image in $\hat
K/U(\hat K)$, 
we may therefore suppose that $K$ is Zariski dense in $H$.

As $H$ is simple, $K$ contains no nontrivial abelian normal
subgroup. Hence $F^*(K)=E(K)$. 

We claim that $K$ contains no nontrivial Zariski closed subgroup
$A$ of $H$. Supposing the contrary, $[E(K),A]$ would be a 
Zariski closed and connected subgroup of $H$ normal in $E(K)$. Such a
group, if nontrivial, would contain an involution, giving a
contradiction. So $[E(K),A]=1$.
 But as $C_{K}(F^*(K))\le Z(F^*(K))$, this shows $A=1$.

In particular, since the characteristic is zero and unipotent subgroups
 are therefore Zariski closed by Lemma \ref{algebraicity_of_unipotence},
 all elements of $K$ are semisimple.
So a Borel subgroup of $K$ is contained in a maximal torus $T$ of $H$.

If $B$ is a Borel subgroup of $H$ containing $T$,
then $\rk(K/T) \leq \rk(H/B)$.
Consider an involution $i\in T$.
Then $H = K C_{H}(i)$ by density of $K$.
So $$ \rk(H) = \rk(C_H(i)) + \rk(K/T) \leq \rk(C_{H}(i)) + \rk(H/B) $$
After cancellation and rearrangement,
 and replacement of $\rk(B)$ by $2+\rk(U)$,
this says $$\dim(U) \leq \dim(C_{H}(i))-2$$
contradicting Table \ref{T:alg_subgrp_dim}.
\end{proof}

\begin{proof}[Proof of Theorem \ref{thm:Altseimer}]
We have $G$ a simple $L^*$-group of finite Morley rank and odd type,
with $M$ strongly embedded.
By Fact \ref{stremb_invconj} and Lemma \ref{Lgrp_oneconjclass} 
we find one of the following
\begin{enumerate}
\item $M$ is a $D$-group;
\item $M^\o/\hat O(M)$ is quasisimple, with central quotient of the
  form $\PSL_2$, $\PSL_3$, or $\G_2$.
\end{enumerate}

In the first case, or in the second case with $M^\o/\hat O(M)$ of the
form $\SL_2$, we have what we claim. We are left therefore with the
following possibilities for the structure of $H = M^\o/\hat O(M)$:
$$\PSL_2, \PSL_3, \SL_3, \G_2$$
By Fact \ref{CK} we have a decomposition of $H$ of the type refuted in Lemma
\ref{CKanalysis} in the cases $H=\PSL_2,\PSL_3,\G_2$, and if $H=\SL_3$
we can pass to $\PSL_3$ by factoring out the center (of odd order).
So in all of these remaining cases we arrive at a contradiction. 
\end{proof}
\section{Strong embedding and $D$-groups}
\label{sec:D-groups}

We aim now at the following. Recall that in the context of
$K^*$-groups of odd type, a simple $D^*$-group is a minimal simple
connected group. 

\begin{theorem}\label{thm:D-group}
Let $G$ be a simple $L^*$-group of finite Morley rank and $p^\perp$ type, 
with $\pr_2(G)\ge 2$.
Suppose that $G$ has a strongly embedded subgroup $M$.
Then the following hold.
\begin{conclusions}
\item $G$ is a $D^*$-group 
\item Sylow 2-subgroups of $G$ are connected.
\end{conclusions}
\end{theorem}

We break the proof into two lemmas.

\begin{lemma}
Let $G$ be a simple $L^*$-group of finite Morley rank and $p^\perp$ type, 
with $\pr_2(G)\ge 2$.
Suppose that $G$ has a strongly embedded subgroup $M$.
Then $M$ is a $D$-group, and Sylow 2-subgroups of $G$ are connected.
\end{lemma}

\begin{proof}
There is a Sylow 2-subgroup $S$ contained in $M$ by Fact \ref{strembdef}.
For any $t\in S$, there is a 2-torus $T$ of $G$ containing $t$
 by Fact \ref{torality}.
So $T \leq M$ by Fact \ref{strembdef}, and $t\in M^\o$.
Thus $S\leq M^\o$ and $\pr_2(M^\o) \geq 2$.
So $M^\o/O(M)$ cannot be $\SL_2$.
Therefore $M$ is a $D$-group by Theorem \ref{thm:Altseimer}.
Now $S$ is connected by Lemma \ref{Sylow_con}.
\end{proof}

\begin{lemma}[{compare \cite[Claim 5.5]{BBN04}}]
\label{Dstar}
Let $G$ be an $L^*$-group of finite Morley rank and odd type.
Suppose that $G$ has a strongly embedded subgroup $M$ which is a $D$-group. 
Then $G$ is a $D^*$-group.
\end{lemma}

\begin{proof}
Suppose that $G$ is a minimal counter example.

We first argue that 
$$ \hbox{$G$ has no proper definable subgroup with an algebraic section.}
  \eqno{(\dag)} $$
So consider a proper definable subgroup $K$ of $G$ with an algebraic section.
As $K$ is not a $D$-group, $K$ contains involutions,
 so we may assume that $M \intersect K$ contains an involution after conjugation.
As $K$ is not contained in $M$, 
 $M \intersect K$ is strongly embedded in $K$, and a $D$-group.
But this contradicts minimality, proving (\dag).
In particular $G$ is connected.

We next show that
$$ \hat O(G) = 1 \mathperiod $$
If $G = M \hat O(G)$, then $G$ is a $D$-group, and so $G > M \hat O(G)$.
We pass to $\bar G = G/\hat O(G)$.
Since $M \hat O(G)$ is strongly embedded in $G$,
 it follows that $\bar M$ is strongly embedded in $\bar G$.
So $\hat O(G) = 1$ follows by minimality.

We now show that 
$$ F^\o(G) = 1 \mathperiod $$
Suppose $F^\o(G) > 1$.  As $\hat O(G)=1$ and $G$ has odd type,
the 2-torsion in $F^\o(G)$ forms a nontrivial 2-torus $T$,
 which is central in $G$.
But $G$ has a strongly embedded subgroup,
 so all involutions are conjugate, and
 hence all involutions of $G$ belong to $T \leq Z(G)$.
It follows easily that $T$ is a Sylow 2-subgroup of $G$.
But then $G$ is a $D$-group after all, a contradiction.

Thus $F(G)$ is finite. and $C(E(G)) = C(F^*(G)) = Z(G)$.
In particular $E(G)$ is nontrivial.
Let $L$ be a quasisimple component of $E(G)$.  Then $Z(L)$ is finite.
So, if $L/Z(L)$ is of degenerate type, then $L$ is of degenerate type
 by Fact \ref{torality}, and $L \leq \hat O(G) = 1$, a contradiction.
So $L/Z(L)$ is algebraic, and hence $L$ is algebraic by \cite{AC99}.
By minimality, we have $G = L$ is quasisimple and algebraic.
But then $G$ contains a subgroup of the form $\SL_2$ or $\PSL_2$,
 and we may suppose this group meets $M$ after conjugation.
As $M$ is a $D$-group this intersection is proper.
So, by minimality, the group $K$ is itself of the form $\SL_2$ or $\PSL_2$.

Now if $G = \SL_2$ then $G \leq C^\o_G(i)$ for $i$ an involution,
 and we contradict strong embedding.
Similarly if $G = \PSL_2$ then
 $M \intersect G < G$ contains a Sylow 2-subgroup of $G$,
but each involution $j \in M$ lies inside
 a torus $T_J \leq C(j) \leq M$ by Fact \ref{torality}.
So $M^\o$ contains $G$, which is again a contradiction.
\end{proof}

\subsection*{Conclusion}
We may now deduce our main results.

\begin{theorem}[$K^*$ Generation Theorem]
\label{K*-generation}
Let $G$ be a simple $K^*$-group of finite Morley rank and $p^\perp$ type, 
with $m_p(G) \geq 3$ and $\pr_2(G) \geq 2$.
Let $V$ be an elementary abelian $p$-group of $p$-rank two. 
Then $\Gamma_V = G$.
In particular
$\Gamma_{S,2}^0 = G$ for any Sylow $p$-subgroup $S$ of $G$.
\end{theorem}

\begin{proof}
If $G$ is not of odd type then it is algebraic of even type by \cite{ABC_EvenType}.
In particular, $G$ is a $K$-group in the pure field language.
It follows that $\Gamma_V = G$ by Theorem \ref{generationA}.
So we may assume that $G$ has odd type too.

Suppose $\Gamma_V<G$.
By Theorems \ref{thm:strongembedding} and \ref{generationA},
 $N(\Gamma_V)$ is strongly embedded in $G$.
By Theorem \ref{thm:D-group}, $G$ is a $D^*$-group,
 hence in this context a minimal connected simple group.
By the main result of \cite{BCJ}, 
a minimal connected simple group of finite Morley rank and
odd type with a proper definable strongly 
embedded subgroup has Pr\"ufer 2-rank one.
This contradiction proves
$$\Gamma_V=G$$

Now by Corollary \ref{cor:2-gencore} we have also $\Gamma_{S,2}^0=G$ for
any Sylow 2-subgroup $S$ of $G$.
\end{proof}

We state the most important case separately.

\begin{corollary}\label{GenerationI}
Let $G$ be a simple $K^*$-group of finite Morley rank and odd type,
with $m_2(G) \geq 3$.
Let $V$ be an elementary abelian 2-group of $2$-rank two. 
Then $\Gamma_V = G$.
In particular $\Gamma_{S,2}^0 = G$ for any Sylow 2-subgroup $S$ of $G$.
\end{corollary}

\begin{proof}
We have $\pr_2(G)\ge 2$ by Theorem \ref{Prank2} and now Theorem
\ref{K*-generation} applies.
\end{proof}

We may also continue the analysis in the $L^*$ case,
 where \cite{BCJ} does not exist.

\begin{theorem} 
\label{L*-generation}
Let $G$ be a simple $L^*$-group of finite Morley rank 
 and odd type,  with $m_2(G) \geq 3$.
Let $V$ be an elementary abelian $2$-group of $2$-rank two with
$\Gamma_V<G$. Then
\begin{enumerate}
\item $G$ is a $D^*$-group.
\item $M := N(\Gamma_V)$ is a strongly embedded subgroup of $G$.
\item The Sylow 2-subgroups of $G$ are connected.
\item $G$ contains a nontrivial unipotent $r$-subgroup
 where $r$ is the least prime divisor of the Weyl group $W$,
 which is nontrivial.
\end{enumerate}
\end{theorem}

For the final point
 we define and use the {\em Weyl group} as in \cite{BC08a}.

\begin{definition}
The {\em Weyl group} $W$ of a group $G$ of finite Morley rank
 is the abstract group $W = N(T)/C\oo(T)$
where $T$ is a maximal decent torus.
\end{definition}

As all maximal decent tori are conjugate by \cite{Ch05},
 the Weyl group is well defined up to conjugacy.

\begin{fact}[{\cite[Theorem 5]{BC08a}}]\label{Weyl_group}
Let $G$ be a connected group of finite Morley rank.
Suppose the Weyl group is nontrivial and has odd order,
 with $r$ the smallest prime divisor of its order.
Then $G$ contains a unipotent $r$-subgroup.
\end{fact}

Here one might prefer $M/M\oo$ over the Weyl group $W$.
Let $T$ denote a Sylow 2-subgroup of $M$.
By a Frattini argument,
 $W_T := N(T)/C\oo(T)$ naturally embeds into $W$.
But $|W| > |W_T|$ seems plausible if $M$ is not solvable.
In particular the minimal prime divisor of $|W|$
 might be strictly less than that of $|M/M\oo|$.
Of course, one may argue that $N_{M\oo }(T)/C\oo_{M\oo}(T) = 1$
 because $M\oo$ is still a $D$-group.
A pair of Frattini arguments then shows that $M/M\oo \cong W_T$.
So one may prove that $M/M\oo$ is nontrivial
 using Fact \ref{stremb_invconj}.

\begin{proof}[Proof of Theorem \ref{L*-generation}]
%
By Theorem \ref{Prank2} we have $\pr_2(G) \ge 2$.
By Theorems \ref{thm:strongembedding} and \ref{generationA},
 $M := N(\Gamma_V)$ is strongly embedded in $G$.
By Theorem \ref{thm:D-group}, $G$ is a $D^*$-group and
 the Sylow 2-subgroups of $G$ are connected.

We now turn our attention to the final point.
Let $T$ be a maximal 2-torus of $M$.
By Fact \ref{stremb_invconj},
 $M$ has only one conjugacy class of involutions.
So $W_T := N(T)/C\oo(T)$ is nontrivial.
By a Frattini argument using \cite{Ch05},
 $W_T$ naturally embeds into $W$.
Hence $W$ is nontrivial as well.
Now $|W|$ has odd order because
 $T$ is a Sylow 2-subgroup of $G$.
So the last point follows from Fact \ref{Weyl_group}.
\end{proof}


\section*{Acknowledgments}

The authors are grateful to Tuna \Altinel, Alexandre Borovik, 
and Ali Nesin for their helpful comments.
We owe special thanks to Christine Altseimer for making her unpublished
work \cite{AltseimerStrEmb} available for use in \S\ref{sec:Altseimer}.
The first author gratefully acknowledges the hospitality and support
of Birmingham, Lyon, Manchester, and Emily Su.
The second author acknowledges the extended and repeated hospitality
of the Universit\'e Claude Bernard (Lyon-I).

\small
\bibliographystyle{alpha} 
\bibliography{burdges,fMr}

\end{document}